\documentclass[12pt]{article}\usepackage{amsmath,amsfonts, amssymb}
 \begin{document}
\newtheorem{proposition}{Proposition}[section]
 \newtheorem{coro}{Corollary}[section]
 \newtheorem{lemma}{Lemma}[section]
 \newtheorem{theorem}{Theorem}[section]
 \newtheorem{definition}{Definition}[section]   
\newfont{\gothic}{eufm10 at 12pt}

\title{ON GEOMETRIC PROPERTIES OF JOINT INVARIANTS OF KILLING TENSORS}

\date{ }                                                             
\author{Caroline M. ADLAM\footnote{Department of Mathematics and Statistics, 
Dalhousie University, Halifax, Nova Scotia, Canada B3H 3J5.}, 
 Raymond G. MCLENAGHAN\footnote{Department of Applied Mathematics, 
University of Waterloo, Waterloo, Ontario, Canada N2L 3G1. The work of the second author
was supported in part by an NSERC Discovery Grant.},\\ and \\
Roman G. SMIRNOV\footnote{Department of Mathematics and Statistics, 
Dalhousie University, Halifax, Nova Scotia, Canada B3H 3J5. The work of the third author
was supported in part by an NSERC Discovery Grant. } }   

\maketitle   
             
\begin{abstract}                   
We employ the language of Cartan's  geometry to present a model for studying vector spaces of Killing two-tensors
defined in pseudo-Riemannian spaces of constant curvature under the action of the corresponding isometry group. 
We also discuss geometric properties of joint invariants of Killing two-tensors defined in the Euclidean plane to formulate and prove
an analogue of the Weyl theorem on joint invariants. In addition, it is shown how the  joint invariants manifest themselves in the theory of superintegrable Hamiltonian systems. 
\end{abstract}

\section{Introduction}\label{section1}
 
The purpose of this article is three-fold. First, we wish to introduce a general model which forms a suitable framework for describing
the invariant theory of Killing tensors from the viewpoint of  Cartan's geometry (see \cite{FO98, FO99,  Gr74, IL03, Sh96} for more details), in particular  its moving frame method and the theory of fiber bundles. Apparently, the idea to study the vector spaces of second order differential operators defined in the Euclidean plane under the action  of the corresponding  isometry group goes back to the classical 1965 Lie group theory paper by Winternitz and  Fri\v{s} \cite{WF65}. In 2002 an analogous idea resurfaced in the problem of studying the  isometry group 
action in a vector space of Killing two-tensors  in the paper by McLenaghan, Smirnov and The \cite{JMP02} on classical Hamiltonian 
systems. A series of recent works  on the subject represent a steady
development of the invariant theory of Killing tensors \cite{CA05, CDM06, HMSD05, hms05, JM05, JMP02, SPT02, MST03, JMP04, MMS04, MSTCUBIC,  SY04, SY05,  RGS, Yue05, YUEPHD}.  

Our second goal is to extend the study of joint invariants of Killing tensors introduced in Smirnov and Yue  \cite{SY04}. More specifically, we formulate and prove an analogue of the Weyl theorem on joint invariants in classical invariant theory (see Olver \cite{PO99} and Weyl \cite{We46} for more details). Furthermore, we introduce the concept of a {\em resultant} of Killing tensors, which is also inspired by its analogue in 
classical invariant theory (see Olver \cite{PO99} for more details). 

Finally, we employ the resultants to derive a new characterization of superintegrable systems. Recall that although examples of superintegrable systems have been known since the time of Kepler, a systematic development of a general  theory of superintegrable systems originated in the pioneering works by Winternitz and collaborators \cite{FMSUW65, MSVW67} (see also Smirnov and Winternitz \cite{SW06} for more references). The theory has been actively developed since then for both quantum and classical Hamiltonian systems. Thus, a general classification and structure theory for superintegrable
systems defined in two- and three-dimensional conformally flat spaces has been introduced in a series of recent papers by 
Kalnins, Kress and Miller (see \cite{KKM4, KKM3, KKM2, KKM1,KKM} and the relevant references therein). As is well-known,  a wide range of  methods are employed to study superintegrable systems. We shall  demonstrate that the geometric methods stemming from the invariant 
theory of Killing tensors also prove their worth in this study and present a new prospective. Thus, the main theorem presented in this work is a geometric characterization of the Kepler potential in terms of the existence of a vanishing resultant of the associated Killing  two-tensors. 

\section{The model}
\label{section2}

In what follows, we shall assume that  $(M, {\bf g})$ is an $m$-dimensional (pseudo-) Riemannian manifold of constant curvature. Generalized Killing tensors can
be defined as the elements of the vector space of solutions to the following overdetermined system of PDEs, called the {\em generalized  Killing 
tensor equation}, given by
\begin{equation}
\label{GKT}
[[\ldots[{\bf K}, {\bf g}], {\bf g}],\ldots, {\bf g}] = 0 \quad \mbox{($n+1$ brackets)},
\end{equation}
where $[$ , $]$ denotes the Schouten bracket \cite{Sch} and ${\bf K}$ is a symmetric contravariant tensor of valence $p$. The generalized 
Killing tensor equation (\ref{GKT}) is determined by given $m \ge 1$, $n \ge 0$ and $p \ge 0$. Let $n$, $m$ and $p$ be fixed and  ${\cal K}_n^p(M)$ denote
the corresponding vector space of solutions to the generalized Killing tensor equation (\ref{GKT}). The {\em Nikitin-Prylypko-Eastwood (NPE)} formula \cite{ME1,ME2,NP90} 
\begin{equation}
\label{NPE}
d = \dim\, {\mathcal K}^p_n(M) = \frac{n+1}{m}{p+ m-1 \choose m-1}{p+n+m    \choose m-1}, 
\end{equation}
represents the dimension $d$ of the vector space ${\cal K}_n^p(M)$. The NPE formula (\ref{NPE}) generalizes the {\em Delong-Takeuchi-Thompson (DTT)}
formula \cite{De, Ta, Th} derived for the case $n=0$. When $(M,{\bf g})$ is not of constant curvature, $d$ given by (\ref{NPE})
is the least upper bound for the dimension of the space of solutions to (\ref{GKT}).
 Of particular importance, due to their relevance in various problems of
mathematical physics, are the elements of the vector space ${\cal K}^2_0(M)\simeq \mathbb{R}^d$ for a given space $(M, {\bf g})$ (see, for example,  \cite{Be97, De, Ei34,  hms05, KKM1, KKM, JMP02, RGS, SW06, WF65} and the relevant references therein). 

The geometric study of   the elements of ${\cal K}^2_0(M)$ having  distinct (and real in the case of $(M,{\bf g})$ being pseudo-Riemannian) eigenvalues and hypersurface forming  eigenvectors (eigenforms) is pivotal in the Hamilton-Jacobi 
theory of orthogonal separation of variables (see \cite{Be97, Ei34, hms05, KKM1, KKM, JMP02, RGS} and the relevant references therein). 
An element ${\bf K} \in {\cal K}^2_0(M)$ enjoying these properties generates an {\em orthogonal coordinate web} (see \cite{Ei34, Be97, hms05} and the relevant references therein) which consists of $m$ foliations, the leaves of which are
$m-1$-dimensional hypersurfaces orthogonal to the eigenvectors of ${\bf K}$. In this way the orthogonal coordinate web is adapted to the hypersurface forming eigenvectors (eigenforms)  of the corresponding Killing tensor $\bf K$. 
The most natural framework for such a study is best described in terms of an appropriate Cartan's (or Klein's) geometry \cite{Sh96, IL03, RGS}, depending on a given (pseudo-) Riemannian manifold $(M,{\bf g})$ and its geometry. 
Indeed, let $G$ be the isometry group of $(M, {\bf g})$, and  $H \subset G$ - a closed subgroup of the Lie group $G$. The study of
the orthogonal coordinate webs generated by the elements of ${\cal K}^2_0(M)$ with the properties prescribed above, is established in the study of the principal fiber $H$-bundle: $\pi_1: \, G \rightarrow G/H  \simeq M$, where the left coset space $G/H$ is identified with 
$M$ on which $G$ acts transitively. We also have the natural structure
of a vector bundle: $\pi_2:\, K^2_0(M) \rightarrow G/H \simeq M.$  The orthogonal coordinate webs are thus defined in the homogeneous space $M$ and the basic equivalence problem can be formulated as follows: Given two orthogonal coordinate webs ${\cal W}_1, {\cal W}_2 \in M$ generated by
the corresponding elements  ${\bf K}_1, {\bf K}_2 \in {\cal K}^2_0(M)$ respectively. The 
coordinate webs ${\cal W}_1, {\cal W}_2$ are said to be {\em equivalent} iff there exists
a group element $g \in G$ and $\ell \in \mathbb{R}$ such that $g{\bf K}_1 = {\bf K}_2 + \ell {\bf g}$. Thus, given two orthogonal coordinate webs ${\cal W}_i$, $i=1,2$, one wishes to know whether or not they are equivalent as defined above.  More specifically, the transitive action $G \circlearrowright M$ induces the corresponding {\em non-transitive} action $G \circlearrowright{\bf K}^2_0(M)$, solving the equivalence problem then amounts to the study of the orbit space $K_0^2(M)/G$. This
structure leads to the principal $G$-bundle $\pi_3: \, {\cal K}^2_0(M) \rightarrow {\cal K}^2_0 (M)/G$, provided the 
action $G \circlearrowright {\bf K}^2_0(M)$ is free (for more details on the fiber bundle theory see, for example,  Fatibene and Francaviglia \cite{FF03}). Finally, one can define
a map $f:\, {\cal K}_0^2(M)/G \rightarrow G$, so that the following diagram commutes.
\begin{center}
\setlength{\unitlength}{1cm}
\begin{picture}(3,3)
\put(0,0){${\cal K}^2_0(M)/G$} 
\put(3,0){${\cal K}^2_0(M)\simeq \mathbb{R}^d$} 
\put(0.7,2){$G$}
\put(3,2){$G/H \simeq M$}
\put(2.9,0.1){\vector(-1,0){1.2}}
\put(0.8,0.4){\vector(0,1){1.4}} 
\put(3.3,0.4){\vector(0,1){1.4}} 
\put(1.1,2){\vector(1,0){1.8}}
\put(0.3,1){$f$}
\put(3.5,1){$\pi_2$}
\put(1.8,2.2){$\pi_1$}
\put(2,0.3){$\pi_3$}
\end{picture} 
\end{center}
\medskip
Now let ${\bf x} \in G/H$ and an element ${\bf K} \in {\cal K}^2_0(M)$ has pointwise distinct eigenvalues and hypersurface forming eigenvectors. If 
${\bf x}$ is not a singular point of the coordinate web (a  point where the eigenvalues of $\bf K$ coincide), we can use the metric $\bf g$ to orthonormalize the eigenvectors $E_1, \ldots, E_m$ of $\bf K$ and then use this {\em frame} of eigenvectors as a basis of $T_{\bf x}G/H$, which has
a group conjugate to $H$ acting on it. Thus, each fiber $\pi_1^{-1}({\bf x})$ can be identified with an orthonormal frame of eigenvectors of $\bf K$. In the other direction the  fiber $\pi_2^{-1}({\bf x})$ can be identified with the corresponding vector space in ${\cal K}^2_0(M)$ for ${\bf x} \in M$. Furthermore, in this vector space we can fix the parameters so that the resulting
element of the vector space is precisely the Killing tensor $\bf K$ above evaluated at ${\bf x} \in M$. It is now a point
in the fiber $\pi^{-1}_2({\bf x})$. Next,  under the projection $\pi_3$ this point is mapped to the corresponding point in the orbit space ${\cal K}_0^2(M)/G$, which is nothing, but the orbit generated by ${\bf  K}$ under the action of $G \circlearrowright {\cal  K}^2_0(M)$. Finally, choosing a function $f$ lifting the action to $G$  is equivalent to either choosing a {\em cross-section} $K$ through the orbits of
${\cal K}_0^2(M)/G$ or fixing a {\em frame}. Then (local) {\em invariants} of the group action $G \circlearrowright {\cal  K}^2_0(M)$
are the coordinates of the canonical forms obtained as the intersection of the orbits with an appropriately chosen
cross-section $K$ \cite{PO99}. In the latter case one can solve  the equivalence problem by employing Cartan's language of the exterior differential forms \cite{Gr74, Sh96, BMS01, IL03}, while in the former case the problem can be solved by using the techniques of the geometric version of the classical  moving frames method as developed by Fels and Olver \cite{FO98, FO99, PO99} (for more details of how these two approaches interact see \cite{RGS}). The equivalence of the two approaches is manifested by the fact that the group $G$ acts transitively on the bundle of frames. We note however, that in practice the most effective approach is to use both methods. Thus, the former method can be used to find canonical forms of the orbits of ${\cal K}_0^2(M)/G$, while the latter approach 
comes into play when, for example, one needs to determine  the moving frames map for  a given representative of an orbit of ${\cal K}_0^2(M)/G$. Note also that the composition map $\gamma = \pi_3\circ f:\, {\cal K}^2_0(M) \rightarrow G$ is the {\em moving frame}
\cite{Gr74, FO98, FO99, PO99}, corresponding to the cross-section (or the frame of eigenvectors) prescribed by a chosen map $f:\, {\cal K}^2_0(M)/G \rightarrow G$. Since 	``The art of doing mathematics consists in finding that special case which contains all the germs of generality.'' (D. Hilbert), in the example that follows we shall do just that.

\section{Illustrative example}
\label{section3} Consider the following example.  $M = \mathbb{E}^2$, $K^2_0(M) = {\cal K}^2_0(\mathbb{E}^2)$, $G = SE(2)$ and $H = SO(2)$. This  example has been extensively studied in the literature from different points of view \cite{BMS01, JM05, JMP02, SY04, RGS, 
WF65}. For the first time the orbit problem $SE(2) \circlearrowright{\cal K}_0^2(\mathbb{E}^2)$ was considered by Winternitz and  Fri\v{s} \cite{WF65} in the context of a classification of quadratic differential operators which commute with the Laplace 
operator. We consider the problem of classification of the orthogonal coordinate webs generated by non-trivial elements of 
the vector space ${\cal K}^2_0(\mathbb{E}^2)$ in the framework of the above diagram. Note first that in this case  a Killing tensor ${\bf K} \in {\cal K}^2_0(\mathbb{E}^2)$ with pointwise distinct eigenvalues has necessarily hypersurface-forming eigenvectors, which generate one of the four orthogonal coordinate webs, namely elliptic-hyperbolic, polar, parabolic and cartesian. The latter depends on 
the corresponding orbit of ${\cal K}^2_0(\mathbb{E}^2)/SE(2)$ that ${\bf K}$ belongs to. More specifically, we fix at ${\bf x} \in \mathbb{E}^2$ the following  orthonormal frame $({\bf x}; E_1, E_2),$
where $E_1,E_2$ are the eigenvectors of ${\bf K}$. We also assume that ${\bf x} \in \mathbb{E}^2$ is not a singular point  of $\bf K$.
   In this view, ${\bf x} = \pi_1:\, SE(2) \rightarrow SE(2)/SO(2)$. The vectors $E_1,E_2$ form a basis of $T_{\bf x}\mathbb{E}^2$. Let $E^1, E^2$
be the corresponding dual basis of $T_{\bf x}^*\mathbb{E}^2$. Then $E^1, E^2$ are
horizontal in the fibration and in this frame the components of  the metric ${\bf g}$ of $\mathbb{E}^2$ and $\bf K$ are given by
\begin{equation}
\label{frame}
g_{ab} = \delta_{ab}E^a\odot E^b \quad \mbox{and} \quad K_{ab} = \lambda_a\delta_{ab}E^a\odot E^b, \quad a,b = 1,2
\end{equation} 
respectively, where $\odot$ is the symmetric tensor product, $\delta_{ab}$ is the Kronecker delta and $\lambda_1$, $\lambda_2$ are the eigenvalues of $\bf K$. Note that
we have used the metric to orthonormalize the frame $(E_1, E_2)$, while  the formulas (\ref{frame}) represent the pull-back under ${\bf x} = \pi_1$ of the metric ${\bf g}$ and the Killing tensor ${\bf K}$ to the bundle of frames. Moreover, we also note that by fixing a frame $(E_1,E_2)$ we have fixed the corresponding map $f:\, {\cal K}_0^2(\mathbb{E}^2)/SE(2) \rightarrow SE(2)$. 
 Then the problem of classification
of the orthogonal coordinate webs generated by the non-trivial elements ${\cal K}^2(\mathbb{E}^2)$ can be solved by making use of the Cartan theory of exterior differential forms \cite{BMS01, RGS}. Thus, upon  introduction of the connection coefficients $\Gamma$  of the Levi-Civita connection $\nabla$ via:
\begin{equation}
\label{nabla}
\nabla_{E_a}E_b = \Gamma_{ab}{}^cE_c, \quad \nabla_{E_c}E^b = \Gamma_{cd}{}^bE^d, \quad a,b, c, d = 1,2,
\end{equation}
we arrive at the  Cartan structure equations
\begin{equation}
\label{MC}
dE^a + \omega^a{}_b \wedge E^b = T^a =0, \quad d\omega^a{}_b + \omega^a{}_c\wedge \omega^c{}_b = \Omega^a{}_b =0, 
\end{equation}
where $\omega^a{}_b : = \Gamma_{cb}{}^aE^c$ are the connection  one-forms, $T^a = \frac{1}{2}T^a{}_{bc}E^b\wedge E^a$
is the vanishing ($\nabla$ is Levi-Civita) torsion two-form and $\Omega^a{}_b $ $:=$ $\frac{1}{2}R^a{}_{bcd}E^c\wedge E^d$ is
the vanishing (in view of flatness of $\mathbb{E}^2$) curvature two-form. The choice of the
connection is not arbitrary. As is well-known from  Riemannian geometry, given a connection $\nabla$ on a manifold $M$ one
can parallel propagate frames. For any path $\tau$ between two points of $M$, parallel transport along $\tau$
defines a linear  mapping $L(\tau)$ between the tangent spaces of two points. This linear map is an {\em isometry} if
the connection $\nabla$ is a Levi-Civita connection. 
Clearly, the linear map $L(\tau)$ induced by a  Levi-Civita connection $\nabla$
maps orthonormal frames to orthonormal frames.  The conditions $E^a\wedge dE^a = 0,$ $ a= 1,2$ are automatically satisfied, hence by (one of the corollaries of) the Frobenius theorem one can introduce curvilinear coordinates $(u,v)$ and functions $f_1$ and $f_2$ such that $E^1 = f_1(u, v)du$ and $E^2 = f_2(u,v)dv$.  Furthermore, it is possible to show that the functions $f_1$ and $f_2$ are such that
$f_1^2 = f_2^2 = A(u) + B(v)$, where $A(u)$ and $B(v)$ are functions of one variable, from which it follows that the metric can be re-written as $ds^2 = (A(u) + B(v))(du^2 + dv^2)$. Note that the explicit forms for $A(u)$ and $B(v)$ 
follow from  the vanishing of the curvature two-form. 
Then the following (differential) {\em invariants} can be used to solve the classification problem (for a complete solution refer to \cite{BMS01,RGS}): 
$\Delta_1(u,v) = -\Gamma_{12}{}^1,$ $\Delta_2(u,v) = -\Gamma_{22}{}^1$. In addition, we introduce in 
this paper the following differential invariant: $\Delta_3(u,v) := \frac{F_{,uu}}{F},$ where $F = \Delta_2/\Delta_1$. The information provided by $\Delta_3(u,v)$  simplifies the classification considerably.  Ultimately, we arrive at the following classification for the orthogonal coordinate webs (orbits)  generated by the non-trivial elements of the vector space ${\cal K}_0^2(\mathbb{E}^2)$.
\begin{equation}
\label{CL1}
\begin{array}{rl}
\Delta_1 = 0, \Delta_2 =0 & \mbox{cartesian}\\
\Delta_1 \Delta_2 =0 & \mbox{polar}\\
\Delta_1\Delta_2 \not=0, \Delta_3=0& \mbox{parabolic}\\
\Delta_1\Delta_2 \not=0,\Delta_3\not= 0 & \mbox{elliptic-hyperbolic}
\end{array}
\end{equation}
However in applications, particularly, those pertinent to the study of Hamiltonian systems of classical mechanics \cite{Be97, Ei34, hms05,JMP02, JMP04, RGS}, the Killing tensors in question arise as functions on the cotangent bundle of the configuration space that define principal parts of first integrals of motion. Commonly the components of such Killing tensors are given in terms of (cartesian) position  coordinates. In this view it is convenient to employ a different version of the moving frames method \cite{FO98, FO99, PO99} upon noticing that the Lie group $SE(2)$ acts transitively on the bundle of frames. We conclude therefore that the same procedure can be carried out ``in the group''. Indeed, let ${\bf x} = (x_1,x_2)$ be cartesian coordinates of $\mathbb{E}^2$. Solving the generalized  Killing tensor equation (\ref{GKT}) for $p=2$, $n = 0$ with respect to these coordinates yields:
\begin{equation}
\begin{array}{rcl}
{\bf K} & = & \displaystyle (\beta_1 + 2\beta_4x_2 + \beta_6x_2^2)\partial_1
\odot  \partial_1 \\ 
& & + \displaystyle (\beta_3 - \beta_4x_1 -
 \beta_5x_2 - \beta_6x_1x_2)\partial_1\odot \partial_2 \\ 
& & + \displaystyle (\beta_2 + 2\beta_5x_1+\beta_6x_1^2) \partial_2\odot
\partial_2,
\end{array}
\label{gKt}
\end{equation}
where $\partial_1 = \frac{\partial}{\partial x_1}$, $\partial_2 = \frac{\partial}{\partial x_2}$ and $\beta_1, \ldots, \beta_6$ are arbitrary constants (of integration) that represent the dimension of the vector space ${\cal K}^2_0(\mathbb{E}^2)$. The three-dimensional Lie group $SE(2)$  of (orientation-preserving) isometries  of $\mathbb{E}^2$ can be represented in matrix form (using the notations adapted  in \cite{IL03}):
$SE(2) = \left\{M\in GL(3,\mathbb{R}) \, \big| \, M=\left(\begin{array}{cc} 1& 0\\ {\bf t} & R \end{array}\right), \, {\bf t} \in \mathbb{R}^2 \, R\in SO(2)\right\}.$
In view of the formula (\ref{gKt}),  the  action $SE(2) \circlearrowright \mathbb{E}^2$ given by 
\begin{equation}
\begin{array}{l}
\tilde{x} = x\cos p_3 - y\sin p_3 + p_1, \\ 
\tilde{y} = x\sin p_3 + y\cos p_3 + p_2 
\end{array}
\label{A1}
\end{equation}
induces the corresponding action  $SE(2) \circlearrowright{\cal K}_0^2(\mathbb{E}^2)$ represented by the following formulas \cite{WF65,JMP02}: 
 \begin{equation}
\label{A2}
\begin{array}{rcl}
\tilde{\beta_1} &=& \beta_1\cos^2p_3 - 2\beta_3\cos p_3\sin p_3 + \beta_2\sin^2p_3 - 2p_2\beta_4\cos p_3 
\\
                  & & - 2p_2\beta_5\sin p_3 + \beta_6p_2^2, \\ 
\tilde{\beta_2} & =& \beta_1\sin^2p_3 - 2\beta_3\cos p_3\sin p_3 +
\beta_2\cos^2p_3 - 2p_1\beta_5\cos p_3  \\ 
 & &+ 2p_1\beta_4\sin p_3  + \beta_6p_1^2,\\
\tilde{\beta_3} & = & (\beta_1-\beta_2)\sin p_3\cos p_3 + \beta_3(\cos^2p_3 - \sin^2p_3) 
\\ 
 & &  + (p_1\beta_4 + p_2\beta_5)\cos p_3 + (p_1\beta_5 - p_2\beta_4)\sin p_3 - \beta_6p_1p_2, \\
 \tilde{\beta_4} & = & \beta_4\cos p_3 + \beta_5\sin p_3 - \beta_6p_2, \\
 \tilde{\beta_5} & = & \beta_5\cos p_3 - \beta_4\sin p_3 - \beta_6p_1,\\
 \tilde{\beta_6} &=& \beta_6. 		  
\end{array}
\end{equation} 
Using the method of moving frames ``\`{a} la Fels and Olver'' \cite{FO98, FO99, PO99}, one arrives at the following (algebraic) {\em invariants}  (see
\cite{RGS} and the relevant references therein for more details):  
\begin{equation}
\label{I23}
\begin{array}{rcl}
{\Delta}'_1 & = & \beta_6,\\
{\Delta}'_2 & = & \beta_6(\beta_1 + \beta_2) - \beta_4^2 - \beta_5^2,\\
{\Delta}'_3 & = & (\beta_6(\beta_1 - \beta_2) - \beta_4^2 + \beta_5^2)^2 + 4(\beta_6 \beta_3 + \beta_4 \beta_5)^2.
\end{array}
\end{equation}
The fact that $\tilde{\Delta}_1$ is an invariant is obvious, the fundamental invariant $\tilde{\Delta}_3$ was derived for the first time in \cite{WF65} (see also \cite{JMP02} where $\tilde{\Delta}_3$ was rederived by making use of Lie's method of infinitesimal generators). Combining the invariants $(\ref{I23})$, one can also make use of the
invariant $k^2$ which is  (half-) the distance between the singular points (foci) in the elliptic-hyperbolic case and is given by 
\begin{equation}
k^2 = \displaystyle \frac{\sqrt{{\Delta}'_3}}{({\Delta}'_1)^2}.
\label{dbf}
\end{equation}
 We also note that the singular points are the foci of the confocal conics (ellipses and hyperbolas) that form the coordinate web in this case. In terms of the algebraic invariants (\ref{I23}) the classification (\ref{CL1}) can be equivalently reformulated as follows  \cite{WF65,JMP02} (see also \cite{SY04}).
\begin{equation}
\label{CL2}
\begin{array}{rl}
{\Delta}'_1 =0,  {\Delta}'_2 =0 & \mbox{cartesian}\\
{\Delta}'_1 = 0, {\Delta}'_3 \not=0 & \mbox{parabolic}\\
{\Delta}'_1 \not=0, {\Delta}'_3=0& \mbox{polar}\\
{\Delta}'_1 \not=0,{\Delta}'_3\not= 0 & \mbox{elliptic-hyperbolic}
\end{array}
\end{equation}

\section{Main theorem}\label{section4}

The next natural step in this study is to consider the action of $G$ on the product
of the vector space and the underlying manifold: ${\cal K}_0^2(M)\times M$  or on $n$ copies of the vector space:  ${\cal K}_0^2(M)\times {\cal K}_0^2(M)\times \cdots \times {\cal K}_0^2(M)$. In the former case the action leads to the {\em covariants} of Killing tensors, while in the latter - the {\em joint invariants} of Killing tensors (see \cite{SY04} for more details). This aspect of the theory also has  much similarity with the study of covariants and joint invariants in classical invariant theory (see Boutin \cite{Bo02},  Olver \cite{PO99}, Weyl \cite{We46} and the relevant  references therein for more details). More specifically, in this work we extend 
the study of joint invariants   introduced in \cite{SY04}. To this end, we shall introduce  an  analogue of the concept of a {\em resultant} in classical invariant theory, as well as to formulate and prove an analogue of the Weyl theorem on joint invariants \cite{PO99, We46} concerning the joint action of $SE(n)$ or $O(n)$ on $n$ copies of $\mathbb{R}^n$: 
$\mathbb{R}^n\times\mathbb{R}^n \times \cdots \times \mathbb{R}^n$. 
\begin{definition} A three-dimensional orbit ${\cal O} \in {\cal K}^2(\mathbb{E}^2)/SE(2)$ is said to be {\em non-degenerate} iff along ${\cal O}$
the invariant $k^2 = \sqrt{{\Delta}'_3}/({\Delta}'_1)^2$ $ \not= 0$, where $\Delta'_1$ and $\Delta'_3$ 
are given by (\ref{CL2}). The action $SE(2) \circlearrowright {\cal K}^2(\mathbb{E}^2)$ for which $k^2 \not=0$ is also said to be {\em non-degenerate}. 
\end{definition}

Recall, that all other orbits of ${\cal K}^2(\mathbb{E}^2)/SE(2)$ are various degeneracies of the non-degenerate orbits defined above \cite{JMP02,RGS}. The (eigenvectors of)  Killing two-tensors corresponding to the non-degenerate orbits generate elliptic-hyperbolic coordinate webs \cite{JMP02,RGS, WF65}. In what follows we shall investigate geometric meaning of the joint invariants of non-degenerate orbits of 
the orbit space $({\cal K}^2(\mathbb{E}^2)\times {\cal K}^2(\mathbb{E}^2) \times \cdots \times {\cal K}^2(\mathbb{E}^2))/SE(2)$. For simplicity we investigate first the joint invariants of the non-degenerate action $SE(2) \circlearrowright {\cal K}^2(\mathbb{E}^2) \times {\cal K}^2(\mathbb{E}^2).$ In this case the group $SE(2)$ acts on each copy of ${\cal K}^2(\mathbb{E}^2)$ with three-dimensional non-degenerate orbits. The group action is given by 
 \begin{equation}
\label{JA}
\begin{array}{rcl}

\tilde{\alpha_1} &=& \alpha_1\cos^2p_3 - 2\alpha_3\cos p_3\sin p_3 + \alpha_2\sin^2p_3 - 2p_2\alpha_4\cos p_3 
\\
                  & & - 2p_2\alpha_5\sin p_3 + \alpha_6p_2^2, \\ 
\tilde{\alpha_2} & =& \alpha_1\sin^2p_3 - 2\alpha_3\cos p_3\sin p_3 +
\alpha_2\cos^2p_3 - 2p_1\alpha_5\cos p_3  \\ 
 & &+ 2p_1\alpha_4\sin p_3  + \alpha_6p_1^2,\\
\tilde{\alpha_3} & = & (\alpha_1-\alpha_2)\sin p_3\cos p_3 + \alpha_3(\cos^2p_3 - \sin^2p_3) 
\\ 
 & &  + (p_1\alpha_4 + p_2\alpha_5)\cos p_3 + (p_1\alpha_5 - p_2\alpha_4)\sin p_3 - \alpha_6p_1p_2, \\
 \tilde{\alpha_4} & = & \alpha_4\cos p_3 + \alpha_5\sin p_3 - \alpha_6p_2, \\
 \tilde{\alpha_5} & = & \alpha_5\cos p_3 - \alpha_4\sin p_3 - \alpha_6p_1,\\
 \tilde{\alpha_6} &=& \alpha_6,\\
\tilde{\beta_1} &=& \beta_1\cos^2p_3 - 2\beta_3\cos p_3\sin p_3 + \beta_2\sin^2p_3 - 2p_2\beta_4\cos p_3 
\\
                  & & - 2p_2\beta_5\sin p_3 + \beta_6p_2^2, \\ 
\tilde{\beta_2} & =& \beta_1\sin^2p_3 - 2\beta_3\cos p_3\sin p_3 +
\beta_2\cos^2p_3 - 2p_1\beta_5\cos p_3  \\ 
 & &+ 2p_1\beta_4\sin p_3  + \beta_6p_1^2,\\
\tilde{\beta_3} & = & (\beta_1-\beta_2)\sin p_3\cos p_3 + \beta_3(\cos^2p_3 - \sin^2p_3) 
\\ 
 & &  + (p_1\beta_4 + p_2\beta_5)\cos p_3 + (p_1\beta_5 - p_2\beta_4)\sin p_3 - \beta_6p_1p_2, \\
 \tilde{\beta_4} & = & \beta_4\cos p_3 + \beta_5\sin p_3 - \beta_6p_2, \\
 \tilde{\beta_5} & = & \beta_5\cos p_3 - \beta_4\sin p_3 - \beta_6p_1,\\
 \tilde{\beta_6} &=& \beta_6, 		  
\end{array}
\end{equation} 
where $\alpha_i,\beta_i,$ $i=1,\ldots, 6$ are the parameters of the respective vector spaces and
the conditions $$k^2_1 = \frac{\sqrt{(\alpha_4^2-\alpha_5^2 + \alpha_6(\alpha_2-\alpha_1))^2 + 4(\alpha_6\alpha_3 + \alpha_4\alpha_5)^2 }}{\alpha_6} \not=0,$$ $$k^2_2 = \frac{\sqrt{(\beta_4^2-\beta5^2 + \beta_6(\beta_2-\beta_1))^2 + 4(\beta_6\beta_3 + \beta_4\beta_5)^2 }}{\beta_6} \not=0$$ hold true. Six joint invariants of the action $SE(2) \circlearrowright {\cal K}^2(\mathbb{E}^2) \times {\cal K}^2(\mathbb{E}^2)$ follow immediately from (\ref{I23}): 

\begin{equation}
\label{I33}
\begin{array}{rcl}
{\Delta}'_1 & = & \alpha_6,\\
{\Delta}'_2 & = & \alpha_6(\alpha_1 + \alpha_2) - \alpha_4^2 - \alpha_5^2,\\
{\Delta}'_3 & = & (\alpha_6(\alpha_1 - \alpha_2) - \alpha_4^2 + \alpha_5^2)^2 + 4(\alpha_6 \alpha_3 + \alpha_4 \alpha_5)^2,\\
{\Delta}'_4 & = & \beta_6,\\
{\Delta}'_5 & = & \beta_6(\beta_1 + \beta_2) - \beta_4^2 - \beta_5^2,\\
{\Delta}'_6 & = & (\beta_6(\beta_1 - \beta_2) - \beta_4^2 + \beta_5^2)^2 + 4(\beta_6 \beta_3 + \beta_4 \beta_5)^2.
\end{array}
\end{equation}
The joint invariants are functionally independent (they are obtained via the moving frames method), any analytic function of the invariants (\ref{I33}) is a joint invariant of the non-degenerate group action $SE(2) \circlearrowright {\cal K}_0^2(\mathbb{E}^2) \times {\cal K}_0^2(\mathbb{E}^2)$ (see Olver \cite{PO99} for more details). In view of the Fundamental Theorem on invariants of regular Lie
group action \cite{PO99}, we need to produce  in total 12 (the dimension of the product space ${\cal K}_0^2(\mathbb{E}^2) \times 
{\cal K}_0^2(\mathbb{E}^2)$) - 3 (the dimension of the orbits) = 9 fundamental joint invariants. To derive the remaining three joint invariants, we employ geometric reasoning. Recall  that each element of a non-degenerate orbit of the action 
$SE(2) \circlearrowright {\cal K}_0^2(\mathbb{E}^2)$ corresponds to a Killing tensor whose eigenvectors generate an elliptic-hyperbolic coordinate web. In turn, the coordinate web is characterized by the foci of the confocal conics (the two families of ellipses and hyperbolas). If $F_1,$
$F_2 \in SE(2)/SO(2) = \mathbb{E}^2$ are the two foci of such a coordinate web, then their respective coordinates $(x_1,y_1)_{F_1}$ and $(x_2,y_2)_{F_2}$ are given in terms of the parameters $\beta_i$, $i = 1,\ldots, 6$ that determine the 
Killing tensor via the formula (\ref{gKt}) by \cite{JMP02}: 
\begin{equation}
\label{F1F2}
\begin{array}{lr}
(x_1,y_1)_{F_1}  = & \\
& \left(\frac{-\beta_5}{\beta_6} + \frac{1}{\beta_6}\left(\frac{\sqrt{\Delta'_6} - \sigma_1}{2}\right)^{1/2}, 
\frac{-\beta_4}{\beta_6} + \frac{1}{\beta_6}\left(\frac{\sqrt{\Delta'_6} + \sigma_1}{2}\right)^{1/2}\right., \\[1cm]
(x_2,y_2)_{F_2}  = & \\
& \left(\frac{-\beta_5}{\beta_6} - \frac{1}{\beta_6}\left(\frac{\sqrt{\Delta'_6} - \sigma_1}{2}\right)^{1/2}, 
\frac{-\beta_4}{\beta_6} - \frac{1}{\beta_6}\left(\frac{\sqrt{\Delta'_6} + \sigma_1}{2}\right)^{1/2}\right., 
\end{array}
\end{equation} 
where $\sigma_1 = \beta_4^2 -\beta_5^2 + \beta_6(\beta_2-\beta_1)$ and $\Delta'_6$ is given by (\ref{I33}). In the case
of the non-degenerate action $SE(2) \circlearrowright {\cal K}_0^2(\mathbb{E}^2)\times {\cal K}_0^2(\mathbb{E}^2)$, we have
two more (generic) foci $F_3$ and $F_4$  corresponding to the second copy of ${\cal K}_0^2(\mathbb{E}^2)$. Their coordinates are
given accordingly by 
\begin{equation}
\label{F3F4}
\begin{array}{lr}
(x_3,y_3)_{F_3}  = & \\
& \left(\frac{-\alpha_5}{\alpha_6} + \frac{1}{\alpha_6}\left(\frac{\sqrt{\Delta'_3} - \sigma_2}{2}\right)^{1/2}, 
\frac{-\alpha_4}{\alpha_6} + \frac{1}{\alpha_6}\left(\frac{\sqrt{\Delta'_3} + \sigma_2}{2}\right)^{1/2}\right., \\[1cm]
(x_4,y_4)_{F_4}  = & \\
& \left(\frac{-\alpha_5}{\alpha_6} - \frac{1}{\alpha_6}\left(\frac{\sqrt{\Delta'_3} - \sigma_2}{2}\right)^{1/2}, 
\frac{-\alpha_4}{\alpha_6} - \frac{1}{\alpha_6}\left(\frac{\sqrt{\Delta'_3} + \sigma_2}{2}\right)^{1/2}\right., 
\end{array}
\end{equation} 
where $\sigma_2 = \alpha_4^2 -\alpha_5^2 + \alpha_6(\alpha_2-\alpha_1)$ and $\Delta'_3$ is given by (\ref{I33}).
Now we can treat the non-degenerate action $SE(2) \circlearrowright {\cal K}_0^2(\mathbb{E}^2)\times {\cal K}_0^2(\mathbb{E}^2)$, 
as the free and regular action $SE(2) \circlearrowright \mathbb{E}^2\times \mathbb{E}^2\times\mathbb{E}^2\times \mathbb{E}^2$, where 
the foci $F_1$, $F_2$, $F_3$ and $F_4$ belong to the corresponding copies of $\mathbb{E}^2$. Recall \cite{PO99, We46}, that
any joint invariant ${\cal J}(F_1,F_2,F_3, F_4)$ of the non-transitive action so defined can be written as a function of the interpoint distances $d(F_1,F_2)$, $d(F_2,F_3)$, $d(F_3,F_4)$, $d(F_4,F_1)$ (this result is also known as ``the Weyl theorem on joint invariants'' \cite{PO99,We46}). Therefore we choose the following square distances as the remaining
three fundamental joint invariants: 
\begin{equation}
\label{I34}
\begin{array}{l}
\Delta'_7 = d^2(F_2,F_3) = (x_2-x_3)^2 + (y_2-y_3)^2, \\
\Delta'_8 = d^2(F_1,F_3) = (x_1-x_3)^2 + (y_1-y_3)^2, \\
\Delta'_9 = d^2(F_2,F_4) = (x_2-x_4)^2 + (y_2-y_4)^2, \\
\end{array}
\end{equation}
where $(x_i,y_i), i = 1,2,3,4$ are specified by (\ref{F1F2}) and (\ref{F3F4}). Direct verification using MAPLE shows that
$\mbox{d}\Delta'_1\wedge\mbox{d}\Delta'_2\wedge \ldots \wedge \mbox{d}\Delta'_9 \not=0$, that is the nine joint invariants given by (\ref{I33})
and $(\ref{I34})$  are functionally independent at a generic point (this conclusion also follows from  a geometric argument). Therefore we arrive at the following result.
\begin{theorem}
Every joint invariant of the non-degenerate action $SE(2) \circlearrowright {\cal K}_0^2(\mathbb{E}^2)\times {\cal K}_0^2(\mathbb{E}^2)$
is a function of the nine fundamental joint invariants $\Delta'_i(\alpha_j,\beta_j)$, $i=1,\ldots,9$, $j = 1,\ldots, 6$ 
given by (\ref{I33}) and (\ref{I34}). 
\label{t31}
\end{theorem} 
Theorem \ref{t31} can naturally be extended to the general case of the non-degenerate group action of $SE(2)$ on $n>2$ copies
of the vector space ${\cal K}_0^2(\mathbb{E}^2)$ as well as, more generally, to the case of Killing two-tensors defined on pseudo-Riemannian spaces of constant curvature of higher dimensions. It must also be mentioned that other joint invariants having a geometric meaning are the various angles, as well as   areas within the quadrilateral $F_1F_2F_3F_4$. For example, the angle $\phi$ given by $$\cos\phi = \frac{\vec{F_1F_3}\cdot
\vec{F_1F_2}}{d(F_1,F_3)d(F_1,F_2)}$$ is a joint invariant. 

We recall now the notion of a {\em resultant} in classical invariant theory (see Olver \cite{PO99} for more details). Thus, let 
\begin{equation} 
\begin{array}{rcl}
P({\bf x}) & = & \tilde{a}_m x^m + \tilde{a}_{m-1}x^{m-1}y + \cdots + \tilde{a}_0y^m, \\
Q({\bf x}) & = & \tilde{b}_n x^n + \tilde{b}_{n-1}x^{n-1}y + \cdots + \tilde{b}_0y^n
\end{array}
\label{PQ}
\end{equation}
be two homogeneous polynomials of degrees $m$ and $n$ respectively. Then a joint invariant of the polynomials given by (\ref{PQ}) is
a function $J(\tilde{\bf a}, \tilde{\bf b})$ of the coefficients $a_0,\ldots, a_m, b_0, \ldots, b_n$ preserved under the action of
$GL(2, \mathbb{R}^2)$ (or its subgroups). Particularly   important joint invariants in this study are those whose vanishing is equivalent to the fact that 
$P$ and $Q$ have common roots. Such joint invariants are said to be {\em resultants} of the system (\ref{PQ}). This important concept can  be transferred naturally to the study of Killing tensors. 
\begin{definition} 
Consider the non-degenerate action $SE(2) \circlearrowright {\cal K}_0^2(\mathbb{E}^2)$ $ \times$ $ {\cal K}_0^2(\mathbb{E}^2)$. Let 
${\bf K}_1,$ ${\bf K}_2$ $\in {\cal K}_0^2(\mathbb{E}^2)$ be two Killing tensors belonging to non-degenerate orbits. Then a {\em 
resultant} ${\cal R}[{\bf K}_1,{\bf K}_2]$ is a joint invariant of the action with the property that the
vanishing ${\cal R}[{\bf K}_1,{\bf K}_2] = 0$ is equivalent to the fact that the orthogonal coordinate webs generated by ${\bf  K}_1$ and ${\bf K}_2$ have a common focus.
\end{definition}

Thus, for example, the joint invariants given by (\ref{I34}) are resultants, while $\Delta'_3 + \Delta'_4$, where $\Delta'_3$ 
and $\Delta'_4$ are given by (\ref{I33}), is not. Another resultant is 
\begin{equation}
\Delta'_{10} = d^2(F_1,F_4) = (x_1-x_4)^2 + (y_1-y_4)^2.
\end{equation} 
In the next section we shall demonstrate how the joint invariants, in particular, the resultants manifest themselves in the study of superintegrable systems. 

\section{An application to the theory of superintegrable systems}
\label{section5}

In order to link the results presented in the preceding sections with the theory of superintegrable systems let us consider now a Hamiltonian
system defined in $\mathbb{E}^2$ by a natural Hamiltonian
\begin{equation} 
\label{Ham}
H({\bf x},{\bf p}) = \frac{1}{2}(p_1^2 + p_2^2) + V({\bf x}), 
\end{equation} 
where ${\bf x} = (x_1,x_2)$ (position coordinates) ${\bf p} = (p_1,p_2)$ (momenta coordinates).
We assume that the Hamiltonian system  (\ref{Ham}) admits two {\em distinct} in a certain sense (see below) first integrals $F_1,F_2$, which are quadratic in the momenta 
\begin{equation}
\label{FI}
F_{\ell}({\bf x},{\bf p}) = K_{\ell}^{ij}({\bf x})p_ip_j + U_{\ell}({\bf x}), \quad {\ell}, i, j = 1,2
\end{equation}
and functionally-independent with $H$. It is well-known that the  vanishing of Poisson brackets $\{F_1, H\} $ $ = $ $\{F_2,H\}$ $  = 0$ yields two sets of conditions, namely the Killing tensor equations
$$[{\bf K}_1, {\bf g} ] = [{\bf K}_2, {\bf g}] = 0,$$
and the compatibility conditions (also known as the Bertrand-Darboux PDEs)
$$\mbox{d}(\hat{\bf K}_1\mbox{d}V) = \mbox{d}(\hat{\bf K}_2\mbox{d}V) =0,$$
where the components of the Killing tensors ${\bf K}_1$, ${\bf K}_2$ are determined by the quadratic in the momenta terms of
the first integrals $F_1$ and $F_2$ respectively given by (\ref{FI}), while the components of the $(1,1)$-tensors  $\hat{\bf K}_1$, $\hat{\bf K}_2$
are as follows: $\hat{K}_{\ell}{}_i^j = \hat{K}_{\ell}{}_{im}g^{mj},$ $\ell = 1,2$, where $g^{ij}$ are the components of the metric
tensor that determines the kinetic part of the Hamiltonian (\ref{Ham}). These assumptions afford orthogonal separation of variables in the associated Hamilton-Jacobi equation of the Hamiltonian system defined by (\ref{Ham}). More specifically, 
the orthogonal coordinates are determined by the eigenvectors (eigenvalues) of the Killing tensors ${\bf K}_1$ and ${\bf K}_2$. Thus, the orthogonal separable coordinate systems can be used to find exact solutions to the Hamiltonian system defined by (\ref{Ham}) via solving the Hamilton-Jacobi equation by separation of variables. Moreover, in view of the above, the Hamiltonian system is superintegrable and multiseparable. 

Furthermore, we assume that the Killing tensors ${\bf K}_1$ and ${\bf K}_2$ belong to two {\em distinct} non-degenerate orbits, namely their respective
eigenvectors (eigenvalues) generate elliptic-hyperbolic coordinate webs. Let $F_1$, $F_2$ be the foci of the first elliptic-hyperbolic coordinate system generated by ${\bf K}_1$, while $F_3$, $F_4$ - the foci of the second one, generated by ${\bf K}_2$. Thus, all four foci are distinct.  We can now treat the pair $\{{\bf K}_1, {\bf K}_2\}$ as an element of the product space ${\cal K}_0^2(\mathbb{E}^2) \times {\cal K}_0^2(\mathbb{E}^2)$. Moreover, it is an element of a non-degenerate orbit of the 
action $SE(2) \circlearrowright {\cal K}_0^2(\mathbb{E}^2) \times {\cal K}_0^2(\mathbb{E}^2)$. Indeed, since $F_1 \not= F_2$ and $F_3 \not = F_4$, the corresponding invariants $k_1^2$ and $k_2^2$ given by (\ref{dbf}) do not vanish. Moreover, $k_1^2 \not= k_2^2$.  In addition, we also have three joint invariants given by (\ref{I34}). In the most general case all of the five joint invariants $k_1^2, k_2^2, \Delta'_7,\Delta'_8, \Delta'_9$  that completely characterize the pair $\{{\bf K}_1, {\bf K}_2\}$ are, in view of Theorem \ref{t31}, functionally independent. At the same time we can treat the Killing tensors ${\bf K}_1$, ${\bf K}_2$, as well as any linear combination $c_1{\bf K}_1 + c_2{\bf K}_2$ of them as elements of the vector space ${\cal K}_0^2(\mathbb{E}^2)$. Consider now the following Killing tensor:
\begin{equation}
\label{generic}
{\bf K}_g = {\bf K}_1 + {\bf K}_2 + \ell {\bf g}, \quad \ell \in \mathbb{R},
\end{equation}
where ${\bf g}$ is the Euclidean metric of the Hamiltonian (\ref{Ham}). Clearly 
\begin{equation}
\mbox{d}(\hat{\bf K}_g\mbox{d}V) =0,
\label{CC}
\end{equation} where
$V$ is the potential part of (\ref{Ham}). Our next observation is that the Killing tensor ${\bf K}_g$ depends upon exactly {\em six} parameters, that is the five joint invariants and $\ell$.  This is only possible if  ${\bf K}_g$ is the most general Killing tensor (\ref{gKt}). But then  the compatibility condition (\ref{CC}) yields that $V$ is constant. 
Since $k_1^2k_2^2 \not=0$, we conclude therefore that one of the joint invariants $\Delta'_7,\Delta'_8, \Delta'_9$
must be a vanishing resultant in order for $V$ to be non-constant. At the same time two of the joint invariants cannot be resultants simultaneously, since $k_1^2 \not= k_2^2$. Without loss of generality let us assume that $\Delta'_8 = {\cal R}[{\bf K}_1, {\bf K}_2] =0$.
It follows that $\Delta'_7 = k^2_1$ and the Killing tensor ${\bf K}_1 + {\bf K}_2$ in (\ref{generic})
depends upon three parameters, namely $k^2_1 = \Delta'_7$, $k^2_1$ and $\Delta'_9$. We also note that since $\Delta'_8$ is a vanishing resultant, the elliptic-hyperbolic coordinate webs generated by ${\bf K}_1$ and ${\bf K}_2$ have a common focus. Now let us set in (\ref{gKt}) the parameters $\beta_3 = 0$ and $\beta_1 = \beta_2$. The resulting Killing tensor is of the form
\begin{equation}
{\bf K}'_g = {\bf K}'_1 + \beta_1{\bf g},
\label{EH2}
\end{equation}
where as before $\bf g$ is the Euclidean metric of the Hamiltonian (\ref{Ham}). Thus, the Killing tensor ${\bf K}'_1$
that appears in (\ref{EH2}) is given by
\begin{equation} 
\begin{array}{rcl}
{\bf K}'_1 & = & \displaystyle (2\beta_4x_2 + \beta_6x_2^2)\partial_1
\odot  \partial_1 \\ 
& & + \displaystyle (- \beta_4x_1 -
 \beta_5x_2 - \beta_6x_1x_2)\partial_1\odot \partial_2 \\ 
& & + \displaystyle ( 2\beta_5x_1+\beta_6x_1^2) \partial_2\odot
\partial_2.
\end{array}
\label{KEH}
\end{equation}
We assume that  the parameters $\beta_4, \beta_5, \beta_6$ are such that both 
$\beta_4^2 + \beta_5^2 \not= 0$ and $\beta_6 \not=0$. 
Note that $\beta_4^2 + \beta_5^2$ and $\beta_6$ are invariants of the action
 $SE(2) \circlearrowright {\cal K}_0^2(\mathbb{E}^2)$ in this case. We also note that in this case the Killing tensor ${\bf K}'_1$ generates an elliptic-hyperbolic web characterized by two foci: $F_1$ with the coordinates $(0,0)$ and $F_2$ with the coordinates $\left(\frac{-2\beta_5}{\beta_6}, \frac{-2\beta_4}{\beta_6}\right)$. Consider now another Killing tensor of the same type, namely given by
 \begin{equation} 
\begin{array}{rcl}
{\bf K}'_2 & = & \displaystyle (2\beta'_4x_2 + \beta'_6x_2^2)\partial_1
\odot  \partial_1 \\ 
& & + \displaystyle (- \beta'_4x_1 -
 \beta'_5x_2 - \beta'_6x_1x_2)\partial_1\odot \partial_2 \\ 
& & + \displaystyle ( 2\beta'_5x_1+\beta'_6x_1^2) \partial_2\odot
\partial_2,
\end{array}
\label{KEH1}
\end{equation}
under the additional assumptions that the three parameters $\beta'_4, \beta'_5, \beta'_6$ are such that $(\beta'_4)^2 + (\beta'_5)^2 \not=0$ and $\beta'_6 \not=0$. Moreover, we assume in addition that $(\beta_4^2 + \beta_5^2)/\beta_6 \not = [( 
\beta'_4)^2 + (\beta'_5)^2]/\beta'_6$, so that the Killing tensors ${\bf K}'_1$, ${\bf K}''_2$ are truly distinct. Note that the coordinates of the  foci
$F_3$, $F_4$ of the elliptic-hyperbolic coordinate web generated by ${\bf K}'_2$ are  $(0,0)$ and $\left(\frac{-2\beta'_5}{\beta'_6}, \frac{-2\beta'_4}{\beta'_6}\right)$ respectively. Thus, it follows that the elliptic-hyperbolic coordinate webs  generated by  ${\bf K}'_1$ and ${\bf K}'_2$ have one common focus at $(0,0)$, or, in other words, the pair $\{{\bf K}'_1, {\bf K}'_2\} \in {\cal K}_0^2(\mathbb{E}^2) \times {\cal K}_0^2(\mathbb{E}^2)$
admits a vanishing resultant.   Therefore  $\{{\bf K}'_1, {\bf K}'_2\}$ have exactly the same geometric properties as  
the pair of Killing tensors $\{{\bf K}_1, {\bf K}_2\}$ given by (\ref{FI}), namely the elliptic-hyperbolic coordinate webs share one focus and the system is characterized by three parameters. Hence, we can {\em identify} ${\bf K}_1 + {\bf K}_2$ with ${\bf K}'_1 + {\bf K}'_2$. Moreover, we see that the sum ${\bf K}'_1 + {\bf K}'_2$ is given by the same formula as either of the Killing tensors (\ref{KEH}) or (\ref{KEH1}). In this view in order to determine the potential ${V}$ in (\ref{Ham}), we substitute the formula (\ref{KEH}) into
the compatibility condition (\ref{CC}) and solve the resulting PDE for $V$. More specifically, after the substitution we set
 successively in (\ref{CC}) $\beta_4 = \beta_5=0$, which
defines a polar web with singular point at $(0,0)$, then $\beta_4=\beta_6=0$, which defines a
parabolic web with singular point at $(0,0)$ and $x$-axis as the focal axis, and finally
$\beta_5=\beta_6=0$, which defines a parabolic web with singular point at $(0,0)$ and
focal axis the $y$-axis.  This results in three PDEs whose unique solution can be easily found. Indeed, the PDEs
enjoy the following forms:
\begin{equation}
2x_2V_{x_1} - 2x_1V_{x_2} - x_1x_2(V_{x_2x_2} - V_{x_1x_1}) + (x_2^2 -x^2_1)V_{x_1x_2} = 0,
\label{PDE1} 
\end{equation}
\begin{equation}
\label{PDE2} 
3V_{x_2}+x_2(V_{x_2x_2} - V_{x_1x_1})+2x_1V_{x_1x_2}=0,   
\end{equation}
\begin{equation}
\label{PDE3}
3V_{x_1}+x_1(V_{x_2x_2} - V_{x_1x_1})+2x_2V_{x_1x_2}=0
\end{equation}
The PDEs (\ref{PDE1})-(\ref{PDE3}) yield 
\begin{equation}
x_2V_{x_1} - x_1V_{x_2} = 0,
\label{PDE4}
\end{equation}
which implies that (in terms of polar coordinates) 
\begin{equation}
\label{PDE5}
V(r, \theta)=V(r\cos\theta,r\sin\theta)
\end{equation}                             
is independent of $\theta$.
Now the general solution of (\ref{PDE3}) in polar coordinates has the form
\begin{equation}
\label{PDE6}
V(r,\theta)=F(r)+\frac{G(\theta)}{r^2}, 
\end{equation}                                    
where $F$ and $G$ are arbitrary functions.  Thus by (\ref{PDE4}) we have
\begin{equation}
\label{PDE7}
G(\theta)=\ell
\end{equation}                                                           
where $\ell$  is an arbitrary constant.  Substituting (\ref{PDE6}) into (\ref{PDE4}) (in polar
coordinates) and taking (\ref{PDE7}) into account, we find that
\begin{equation}
F(r)=-\frac{{\ell}}{r^{2}}-\frac{m }{r} + n,                                           
\label{PDE8}
\end{equation}
where $\ell$ as in (\ref{PDE7}) and $m,n$ are arbitrary constants. Substituting (\ref{PDE7}) and (\ref{PDE8}) into
(\ref{PDE6}) and transforming back to cartesian coordinates, we arrive at the potential of the {\em Kepler
problem}:
\begin{equation}
\label{K}
V(x_1,x_2) = \frac{1}{\sqrt{x_1^2 + x_2^2}},
\end{equation}
where without loss of generality we set $m = -1$, $n=0$. 
The inverse problem, namely when one starts with the Kepler potential (\ref{K})  and then finds the most general
Killing tensor compatible with (\ref{K}) via (\ref{CC}), can be solved in the same manner, that is by solving the corresponding (Bertrand-Darboux) PDE. The calculations are straightforward and we present here the result only. 
Thus, the most general Killing tensor compatible with the Kepler potential (\ref{K}) is precisely the three-parameter family of Killing tensors 
given by (\ref{KEH}). We conclude therefore that we have proven the following theorem. 
\begin{theorem}
\label{t52}
Let the potential $V$ of the general Hamiltonian (\ref{Ham}) be compatible via (\ref{CC}) with any two non-degenerate
Killing tensors ${\bf K}_1, {\bf K}_2 \in {\cal K}_0^2(\mathbb{E}^2)$. Then the following statements are equivalent. 
\begin{itemize}
\item[(1)] The pair of Killing tensors $\{{\bf K}_1,$ ${\bf K}_2 \} \in {\cal K}_0^2(\mathbb{E}^2)\times {\cal K}_0^2(\mathbb{E}^2)$ admits one vanishing resultant ${\cal R}[{\bf K}_1, {\bf K}_2]$.

\item[(2)] The potential $V$ given by (\ref{Ham}) is the Kepler potential (\ref{K}). 
\end{itemize}
\end{theorem}

\section{Conclusions}\label{section6}
The results presented in this article can naturally be generalized by increasing the dimension of the underlying (pseudo-) Riemannian manifold, changing its curvature and the signature of its metric. Some of these cases will be investigated by the authors in the forthcoming papers on the subject.

\end{document}